\documentclass[12pt,twoside,reqno]{amsart}
\usepackage{amsxtra,amscd}
\usepackage{graphicx}
\usepackage{amsmath}
\usepackage{amsfonts}
\usepackage{amssymb}

\newtheorem{theorem}{Theorem}[section]
\newtheorem{lemma}[theorem]{Lemma}
\theoremstyle{definition}
\newtheorem{definition}[theorem]{Definition}

\newtheorem{remark}[theorem]{Remark}

\numberwithin{equation}{section}

\begin{document}
\title{on the existence of geometric models for function fields in
several variables}
\author{Feng-Wen An}
\address{School of Mathematics and Statistics, Wuhan University, Wuhan,
Hubei 430072, People's Republic of China}
\email{fwan@amss.ac.cn}

\begin{abstract}
In this paper we will give an explicit construction of the geometric model for a prescribed Galois extension of a function field in several variables over a number field. As  a by-product,
we will also prove the existence of quasi-galois closed covers of arithmetic schemes, which is a generalization of the pseudo-galois covers of arithmetic varieties in the sense of Suslin-Voevodsky.
\end{abstract}

\maketitle

\begin{center}
{\tiny {Contents} }
\end{center}

{\tiny \qquad {Introduction} }

{\tiny \qquad {1. Statement of The Main Theorem} }

{\tiny \qquad {2. An Explicit Construction for the Model} }

{\tiny \qquad {3. Proof of The Main Theorem} }

{\tiny \qquad {References}}

\section*{Introduction}

Let $F$ be a number field and let $E$ be a finitely generated extension of $F$.  If $tr.deg_{F}E=1$, this is on the theory of function fields of one variable, especially on the Riemann-Roch Theory.

Consider the case that $[E:F]<\infty$. In recent decades one
 has been attempted to use the related data of arithmetic varieties $X/Y$ to
describe such a given (Galois) extension $E/F$ (for example, see \cite{Bloch,Kato,Schmidt,Lang,Raskind,Saito,VS1,w1,w2}). The reason is that there is a nice relationship between them:

For the case that $\dim X=\dim Y$, as it has been seen, under certain conditions the
arithmetic varieties $X/Y$ behave like Galois extensions $E/F$ of number fields; at the same time,
their automorphism groups $Aut(X/Y)$ behave like the Galois groups $Gal(E/F)$. In particular, the related data of varieties, such as the arithmetic fundamental groups, encode plenty of information of the maximal abelian class fields of the number fields. It needs to decode them for one to obtain class fields.

Moreover, in such a case, one says that the arithmetic varieties $X/Y$ are a \emph{geometric model} for the Galois extensions $E/F$ if the Galois group $Gal\left( E/F\right) $
is isomorphic to the automorphism group ${Aut}\left( X/Y\right) $ (for
example, see \cite{SGA1,GIT,Raskind,VS1,SV2}).

Now let $F$ be a finitely generated extension over a number field.
In this paper we will have a try to use the related data of arithmetic varieties $X/Y$ to
describe a prescribed field $E$, a finitely generated extension over $F$, of transcendental degree not less than one. We will give an explicit construction of such a geometric model for function fields in several
variables (see \emph{Main Theorem}).

On the other hand, we will also demonstrate the existence of quasi-galois closed covers in \cite{An2}, which is as a by-product of the procedure for the proof of the \emph{Main Theorem} (see \emph{Theorem 3.4}). It can be regarded as a generalization of the pseudo-galois covers of arithmetic varieties in the sense of Suslin-Voevodsky (see \cite{VS1,SV2}).

\textbf{Acknowledgment.} The author would like to express his sincere gratitude to Professor Li
Banghe for his invaluable advice and instructions on algebraic geometry and
topology.

\section{Statement of the Main Theorem}

In the present paper, an \textbf{arithmetic variety} is an integral scheme of
finite type over $Spec\left( \mathbb{Z}\right) $. Let $k(X)\triangleq \mathcal{O}_{X,\xi }$ denote the function field of an arithmetic variety $X$
with generic point $\xi $.

Let $E$ be a finitely generated extension of a field $F$. Here $E$ is not necessarily algebraic over $F$. Then $E $ is said to be a \textbf{Galois extension} of $F$ if $F$ is the fixed subfield of the Galois group $Gal\left( E/F\right) $ in $E$.

The following is the main theorem of the paper.

\begin{theorem}
\emph{(\textbf{The Main Theorem})} Let $K$ be a finitely generated extensions over a number
field. Suppose that $Y$ is an arithmetic variety with $K=k\left( Y\right)$.
Take any finitely generated extensions $L$ of $K$ such that $L$ is
Galois over $K$.

Then there exists an arithmetic variety $X$ and a surjective morphism $
f:X\rightarrow Y$ of finite type such that
\begin{itemize}
\item $L=k\left( X\right) $;
\item the morphism $f$ is affine;
\item there
is a group isomorphism $Aut\left( X/Y\right) \cong Gal\left( L/K\right) $.
\end{itemize}
\end{theorem}

\begin{remark}
Let $\dim X=\dim Y$. Then $X/Y$ are said to be
 a \textbf{geometric model} of the field extension $E/F$ provided that $k\left( X\right) =E$ and $%
k\left( Y\right) =F$ and there is a group isomorphism $Aut\left( X/Y\right)
\cong Gal\left( E/F\right) $ (for example, see  \cite{GIT,Raskind,VS1,SV2}). In the paper \emph{Theorem 1.1} above gives us a geometric model for function fields in several variables, which is an analogue of the case for finite extensions of number fields.
\end{remark}

\begin{remark}
In the course of the proof of  \emph{Theorem 1.1}, as  a by-product,
we will also demonstrate the existence of quasi-galois closed covers in \cite{An2} (see \emph{Theorem 3.4} in the paper), which can be regarded as a generalization of the pseudo-galois covers of arithmetic varieties in the sense of Suslin-Voevodsky (see \cite{VS1,SV2}).

\end{remark}

\section{An Explicit Construction for the Model}

\subsection{Notation}

Let us fix some notation and definitions before we give the procedure for the
construction (for details, see \cite{An1,An2}). Given an integral domain $D$. Let $Fr(D)$ denote the field of fractions on $D$. In particular, if $D$ is a subring of
a field $\Omega $, the field $Fr(D)$ will always assumed to be contained in $\Omega $.

Let $(X,\mathcal{O}_{X})$ be a scheme. As usual, an affine covering of the
scheme $(X,\mathcal{O}_{X})$ is a family $\mathcal{C}_{X}=\{(U_{\alpha
},\phi _{\alpha };A_{\alpha })\}_{\alpha \in \Delta }$ such that for each $%
\alpha \in \Delta $, $\phi _{\alpha }$ is an isomorphism from an open set $%
U_{\alpha }$ of $X$ onto the spectrum $Spec{A_{\alpha }}$ of a commutative
ring $A_{\alpha }$. Each $(U_{\alpha },\phi _{\alpha };A_{\alpha })\in
\mathcal{C}_{X}$ is called a \textbf{local chart}. An affine covering $%
\mathcal{C}_{X}$ of $(X,\mathcal{O}_{X})$ is said to be \textbf{reduced} if $%
U_{\alpha }\neq U_{\beta }$ holds for any $\alpha \neq \beta $ in $\Delta $.

Let $\mathfrak{Comm}$ be the category of commutative rings with identity.
Fixed a subcategory $\mathfrak{Comm}_{0}$ of $\mathfrak{Comm}$. An affine
covering $\{(U_{\alpha},\phi_{\alpha};A_{\alpha })\}_{\alpha \in \Delta}$ of
$(X, \mathcal{O}_{X})$ is said to be \textbf{with values} in $\mathfrak{Comm}%
_{0}$ if
 for each $\alpha \in \Delta $ there are $\mathcal{O}_{X}(U_{\alpha})=A_{\alpha}$ and $U_{\alpha}=Spec(A_{\alpha})$, where
 $A_{\alpha }$ is a ring contained in $\mathfrak{Comm}_{0}$.

Let $\Omega $ be a field and let $\mathfrak{Comm}(\Omega )$ be the category consisting of the subrings of $\Omega $ and
their isomorphisms. An affine covering $\mathcal{C}_{X}$ of $(X,\mathcal{O}%
_{X})$ with values in $\mathfrak{Comm}(\Omega )$ is said to be \textbf{with
values in the field $\Omega $}.

\subsection{Process of the Construction}

The following is the procedure for the construction of the geometric model.

Let $K$ be a finitely generated extensions over a number field and let $Y$
be an arithmetic variety such that $K=k\left( Y\right) .$ Take any finitely
generated extensions $L$ of $K$ such that $L/K$ is a Galois extension.

We will proceed in several steps to construct an arithmetic
variety $X$ and a surjective morphism $f:X\rightarrow Y$ satisfying the desired property in the \emph{Main Theorem} of the paper, which will be proved in next section.

\emph{\textbf{Step 1.}} Fixed an algebraic closure $\Omega _{L}$ of $L.$ Put $$\Omega
_{K}=\Omega _{L}\cap \overline{K},$$ i.e., an algebraic closure of $K.$

Without loss of generality, assume that the ring $\mathcal{O}_{Y}\left(
V\right) $ is contained in $\Omega _{K}$ for each affine open set $V$ of the
scheme $Y.$

Otherwise, if that property does not hold, by discussion in \cite{An1} we can choose a
scheme $\left( Y^{\prime },\mathcal{O}_{Y^{\prime }}\right) $ which has that
property and is isomorphic to $\left( Y,\mathcal{O}_{Y}\right) $.

Evidently, that property holds automatically if $Y$ is an affine scheme.

Choose the elements
$$
t_{1},t_{2},\cdots ,t_{n}\in L\setminus K
$$
to be a nice basis of $L$ over $K$ (see \cite{An2}), that is, they satisfy the
following conditions:

$\left( i\right) $ $L=K(t_{1},t_{2},\cdots ,t_{n})$;

$\left( ii\right) $ $t_{1},t_{2},\cdots ,t_{r}$ constitute a transcendental
basis of $L$ over $K$;

$\left( iii\right) $ $t_{r+1},t_{r+2},\cdots ,t_{n}$ are linearly
independent over $K(w_{1},w_{2},\cdots ,w_{r})$, where $0\leq r\leq n$.

Let $\mathcal{C}_{Y}$ be the maximal element (by set inclusion) in the
collection of the reduced affine coverings of the scheme $Y$ with values in $\Omega _{K}.$

\emph{\textbf{Step 2.}} Take any local chart $\left(
V,\psi _{V},B_{V}\right) \in \mathcal{C}_{Y}.$ Then $V$ is an affine open subset of $Y $ and we have $$Fr\left(
B_{V}\right) =K\text{ and } \mathcal{O}_{Y}\left( V\right) =B_{V}\subseteq \Omega
_{K}.$$

Define $A_{V}$ to be the subring of $L$ generated over $B_{V}$ by the
set of elements in $L$
$$\Delta _{V}\triangleq \{\sigma \left( t_{j}\right) \in L:\sigma \in
Gal\left( L/K\right) ,1\leq j\leq n\}.$$

That is, we have $$A_{V}=B_{V}\left[ \Delta _{V}\right] .$$

Put $$\Delta _{V}^{\prime
}=\Delta _{V}\setminus \{t_{1},t_{2},\cdots ,t_{r}\}.$$ We have
$$Fr\left( A_{V}\right) =L;$$
$$A_{V}=B_{V}\left[ t_{1},t_{2},\cdots ,t_{r}\right] \left[ \Delta
_{V}^{\prime }\right] .$$

Then $\Delta _{V}^{\prime }$ is
a nonvoid set. It is seen that $B_{V}$ is exactly the invariant subring of
the natural action of the Galois group $Gal\left( L/K\right) $ on $A_{V}.$

Set $$i_{V}:B_{V}\rightarrow A_{V}$$ to be the inclusion.

\emph{\textbf{Step 3.}} Define the disjoint union $$\Sigma =\coprod\limits_{\left(
V,\psi _{V},B_{V}\right) \in \mathcal{C}_{Y}}Spec\left( A_{V}\right) .$$

Let $$\pi _{Y}:\Sigma \rightarrow Y$$ be the projection.

$\Sigma $ is a topological space, where the topology $\tau _{\Sigma }$ on $\Sigma $ is naturally
determined by the Zariski topologies on all $Spec\left( A_{V}\right) .$

\emph{\textbf{Step 4.}} Define an equivalence relation $R_{\Sigma}$ in $\Sigma$ in such a manner:

Take any $x_{1},x_{2}\in \Sigma $. We say $$x_{1}\sim x_{2}$$ if and only if $$
j_{x_{1}}=j_{x_{2}}$$ holds in $L$.

Here, $j_{x}$ denotes the corresponding prime
ideal of $A_{V}$ to a point $x\in Spec\left( A_{V}\right) $ (see \cite{EGA}).

Define $$X=\Sigma /\sim .$$ Let $$\pi _{X}:\Sigma \rightarrow X$$ be the
projection.

Hence, $X$ is a
topological space as a quotient of $\Sigma .$

\emph{\textbf{Step 5.}} Define a map $$f:X\rightarrow Y$$ by $$\pi _{X}\left( z\right)
\longmapsto \pi _{Y}\left( z\right) $$ for each $z\in \Sigma $.

\emph{\textbf{Step 6.}} Put $$\mathcal{C}
_{X}=\{\left( U_{V},\varphi _{V},A_{V}\right) \}_{\left( V,\psi
_{V},B_{V}\right) \in \mathcal{C}_{Y}}$$ where
  $U_{V}=\pi _{Y}^{-1}\left( V\right) $ holds and $\varphi _{V}:U_{V}\rightarrow Spec(A_{V})$ is the identity map for each $\left(
V,\psi _{V},B_{V}\right) \in \mathcal{C}_{Y}$.
Then $\mathcal{C}
_{X}$ is a reduced affine covering on the
space $X$ with values in $\Omega _{L}$.

Define the scheme $$\left( X,\mathcal{O}_{X}\right) $$ to be  obtained  by gluing the affine schemes $Spec\left(
A_{V}\right) $ for all local charts $\left( V,\psi _{V},B_{V}\right) \in \mathcal{C}_{Y}$
with respect to the equivalence relation $R_{\Sigma}$ (see \cite{EGA,Hrtsh}).

Then $\mathcal{C}_{X}$ is admissible and the sheaf $\mathcal{O}_{X}$ is
an extension of $\mathcal{C}_{X}$ on the space $X$ (see \cite{An1}).

Finally, $\left( X,\mathcal{O}_{X}\right) $ is the desired scheme and
 $f:X\rightarrow Y$ is the desired morphism of schemes. (Note that the proof will be given in the following section.)

This completes the construction.

\section{Proof of the Main Theorem}

\subsection{Definitions}

Assume that $\mathcal{O}_{X}$ and $\mathcal{O}^{\prime}_{X}$ are two structure sheaves on the underlying space of an integral scheme $X$. The integral schemes $(X,\mathcal{O}_{X})$ and $(X, \mathcal{O}^{\prime}_{X})$ are said to be \textbf{essentially equal} provided that for any open set $U$ in $X$, we have
 $$U \text{ is affine open in }(X,\mathcal{O}_{X}) \Longleftrightarrow \text{ so is }U \text{ in }(X,\mathcal{O}^{\prime}_{X})$$ and in such a case,  $D_{1}=D_{2}$ holds or  there is $Fr(D_{1})=Fr(D_{2})$ such that for any nonzero $x\in Fr(D_{1})$, either $$x\in D_{1}\bigcap D_{2}$$ or $$x\in D_{1}\setminus D_{2} \Longleftrightarrow x^{-1}\in D_{2}\setminus D_{1}$$ holds, where $D_{1}=\mathcal{O}_{X} (U)$ and $D_{2}=\mathcal{O}^{\prime}_{X} (U)$.

 Two schemes $(X,\mathcal{O}_{X})$ and $(Z,\mathcal{O}_{Z})$ are said to be \textbf{essentially equal} if the underlying spaces of $X$ and $Z$ are equal and the schemes $(X,\mathcal{O}_{X})$ and $(X,\mathcal{O}_{Z})$ are essentially equal.

Let $X$ and $Y$ be two arithmetic varieties and let $f:X\rightarrow Y$ be a
surjective morphism of finite type. By a \textbf{conjugate} $Z$ of $X$ over $
Y$ we understand an arithmetic variety $Z$ that is isomorphic to $X$ over $Y$
. Let $Aut\left( X/Y\right) $ denote the group of automorphisms of $X$ over $
Z$.

Then
$X$ is said to be \textbf{quasi-galois closed} over $Y$ by $f$ if there is an algebraically closed field $\Omega$
and a reduced affine covering $\mathcal{C}_{X}$ of $X$ with values in $
\Omega $ such that for any conjugate $Z$ of
$X$ over $Y$ the two conditions are satisfied:
\begin{itemize}
\item $(X,\mathcal{O}_{X})$ and $(Z,\mathcal{O}_{Z})$ are essentially equal if $Z$ has a reduced
affine covering with values in $\Omega$.

\item $\mathcal{C}_{Z}\subseteq \mathcal{C}_{X}$ holds if $\mathcal{C}_{Z}$
is a reduced affine covering of $Z$ with values in $\Omega $.
\end{itemize}

Let $K$ be an extension of a field $k$. Here $K/k$ is not necessarily
algebraic. Recall that $K$ is said to be
\textbf{quasi-galois} over $k$ if each irreducible polynomial $f(X)\in F[X]$
that has a root in $K$ factors completely in $K\left[ X\right] $ into linear
factors for any intermediate field $k\subseteq F\subseteq K$ (see \cite{An2}).

The elements
$$
t_{1},t_{2},\cdots ,t_{n}\in K\setminus k
$$
to be a \textbf{nice basis} of $K$ over $k$  if they satisfy the
following conditions:

$\left( i\right) $ $L=K(t_{1},t_{2},\cdots ,t_{n})$;

$\left( ii\right) $ $t_{1},t_{2},\cdots ,t_{r}$ constitute a transcendental
basis of $L$ over $K$;

$\left( iii\right) $ $t_{r+1},t_{r+2},\cdots ,t_{n}$ are linearly
independent over $K(t_{1},t_{2},\cdots ,t_{r})$, where $0\leq r\leq n$.

Now let $D\subseteq D_{1}\cap D_{2}$ be three integral domains.
The ring $D_{1}$ is said to be \textbf{quasi-galois} over $D$ if the field $Fr\left( D_{1}\right) $ is a
quasi-galois extension of $Fr\left( D\right) $.

The ring $D_{1}$ is said to be a  \textbf{conjugation} of $D_{2}$ over $D$ if there is a $(r,n)-$nice $k-$basis
$w_{1},w_{2},\cdots ,w_{n}$ of the field $Fr(D_{1})$ and an
$F-$isomorphism $\tau_{(r,n)}:Fr(D_{1})\rightarrow Fr(D_{2})$ of
fields such that  $$\tau_{(r,n)}(D_{1})=D_{2},$$ where $k=Fr(D)$ and
 $F\triangleq k(w_{1},w_{2},\cdots ,w_{r})$ is assumed to be contained in
the intersection $Fr(D_{1})\cap Fr(D_{2})$.

\subsection{Criterion for Quasi-gaois Closed}

Let $X$
and $Y$ be two arithmetic varieties. Let $\Omega $ be a fixed algebraically closed closure of the function
field $k\left( X\right) $.

\begin{definition}
Let $\varphi :X\rightarrow Y$ be a surjective
morphism of finite type. A reduced affine covering
$\mathcal{C}_{X}$ of $X$ with values in $\Omega $ is said to be \textbf{%
quasi-galois closed} over $Y$ by $\varphi$ if the below condition is satisfied:

There exists a local chart $(U_{\alpha }^{\prime },\phi
_{\alpha }^{\prime };A_{\alpha }^{\prime })\in \mathcal{C}_{X}$ such that $U_{\alpha }^{\prime }\subseteq \varphi^{-1}(V_{\alpha})$ for any
 $(U_{\alpha },\phi _{\alpha };A_{\alpha })\in \mathcal{C}_{X}$, for any affine open set $V_{\alpha}$ in $Y$ with $U_{\alpha }\subseteq \varphi^{-1}(V_{\alpha})$,
and for any conjugate $A_{\alpha }^{\prime }$ of $A_{\alpha }$ over $B_{\alpha}$, where $B_{\alpha}$ is the canonical image of $\mathcal{O}_{ Y}(V_{\alpha})$ in the function field $k(Y)$.
\end{definition}

\begin{lemma}
Let $\varphi :X\rightarrow Y$ be a surjective morphism of finite type. Suppose that the function field $k(Y)$ is contained in $\Omega$.
Then the scheme $X$ is quasi-galois closed over $Y$ if there is a unique maximal reduced affine covering $\mathcal{C}_{X}$ of $X$ with values in $\Omega $ such that $\mathcal{C}_{X}$ is
quasi-galois closed over $Y.$
\end{lemma}

\begin{proof}
Assume that there is a unique maximal reduced affine covering $\mathcal{C}_{X}$ of $X$ with values in $\Omega $ such that $\mathcal{C}_{X}$ is
quasi-galois closed over $Y.$

Fixed any a conjugate $Z$ of $X$ over $Y$. Let $\sigma:Z \rightarrow X$ be an isomorphism of schemes over $Y$. Suppose that $Z$ has a reduced affine covering $\mathcal{C}_{Z}$ with values in $\Omega$.

Take any local chart $(W,\delta,C)\in \mathcal{C}_{Z}$. Put $$U=\sigma (W);$$
$$A=\mathcal{O}_{X}(U);$$  $$C=\mathcal{O}_{Z}(W).$$
Then we have $$U=Spec(A) \text{ and } W=Spec(C).$$

As $\mathcal{C}_{X}$ is quasi-galois closed over $Y$, it is seen that there is an affine open subset $U^{\prime}$ in $X$ such that $$C=\mathcal{O}_{X}(U^{\prime}).$$ As $$U^{\prime}=Spec(C)=W,$$ we have $$\sigma^{-1}(U)=U^{\prime}\subseteq X;$$ hence, $$Z=\sigma^{-1}(X)=X.$$ It follows that we must have $(X,\mathcal{O}_{X})=(Z,\mathcal{O}_{Z}).$
\end{proof}

An affine covering $\{(U_{\alpha },\phi _{\alpha };A_{\alpha })\}_{\alpha
\in \Delta }$ of $(X,\mathcal{O}_{X})$ is said to be an \textbf{affine
patching} of $(X,\mathcal{O}_{X})$ if $\phi _{\alpha }$ is the identity map
on $U_{\alpha }=SpecA_{\alpha }$ for each $\alpha \in \Delta .$

Evidently,
an affine patching is reduced.

\begin{lemma}
Let $\varphi :X\rightarrow Y$ be a
surjective morphism of finite type. Suppose that the function field $k(Y)$ is contained in $\Omega$. Then $X$ is quasi-galois closed over $Y$ if there is a unique maximal
affine patching $\mathcal{C}_{X}$ of $X$ with values in $\Omega $ such that
\begin{itemize}
\item either $\mathcal{C}_{X}$ is
quasi-galois closed over $Y$,

\item or $A_{\alpha }$ has only one conjugate over $B_{\alpha}$ for any
 $(U_{\alpha },\phi _{\alpha };A_{\alpha })\in \mathcal{C}_{X}$ and for any affine open set $V_{\alpha}$ in $Y$ with $U_{\alpha }\subseteq \varphi^{-1}(V_{\alpha})$, where $B_{\alpha}$ is the canonical image of $\mathcal{O}_{ Y}(V_{\alpha})$ in the function field $k(Y)$.
\end{itemize}
\end{lemma}

\begin{proof}
It is immediate from \emph{Lemma 3.2}.
\end{proof}

\subsection{Existence of Quasi-gaois Closed Covers} Now we give the existence of quasi-galois closed covers which take values in a prescribed extension of the function field in several variables.

\begin{theorem}
Let $K$
be a finitely generated extensions of a number field and let $Y$ be an
arithmetic variety with $K=k\left( Y\right)$. Fixed any finitely
generated extensions $L$ of $K$ such that $L$ is Galois over $K$.

Then there exists an arithmetic variety $X$ and a surjective morphism $
f:X\rightarrow Y$ of finite type such that
\begin{itemize}
\item $L=k\left( X\right) $;
\item the morphism $f$ is affine;
\item $X$ is
a quasi-galois closed over $Y$ by $f$.
\end{itemize}
\end{theorem}

\begin{proof}
It is immediate from \emph{Lemma 3.3} and the construction in \S 2.
\end{proof}

\subsection{Proof of the Main Theorem}

Now we can give the proof of the Main Theorem of the paper.

\begin{proof}
\textbf{(Proof of Theorem 1.1)} It is immediate from \emph{Theorem 3.4} above and the \emph{Main Theorem}
in \cite{An2}.
\end{proof}


\newpage

\begin{thebibliography}{99}
\bibitem{An1} An, F-W. The Affine Structures on a Ringed Space and Schemes.
eprint arXiv:0706.0579.

\bibitem{An2} An, F-W. Automorphism groups of quasi-galois closed arithmetic
schemes. eprint arXiv:0907.0842.

\bibitem{Bloch} Bloch, S. Algebraic $K-$Theory and Classfield Theory for Arithmetic Surfaces. Annals of Math, 2nd Ser., Vol 114,
No. 2 (1981), 229-265.

\bibitem{EGA} Grothendieck, A; Dieudonn\'{e}, J. \'{E}l\'{e}ments de G\'{e}%
oem\'{e}trie Alg\'{e}brique. vols I-IV, Pub. Math. de l'IHES, 1960-1967.

\bibitem{SGA1} Grothendieck, A; Raynaud, M. Rev$\hat{e}$tements $\acute{E}$%
tales et Groupe Fondamental (SGA1). Springer, New York, 1971.

\bibitem{Hrtsh} Hartshorne, R. Algebraic Geometry. Springer, New York, 1977.

\bibitem{Kato} Kato, K; Saito, S. Unramified Class Field
Theory of Arithmetical Surfaces. Annals of Math, 2nd Ser., Vol 118,
No. 2 (1983), 241-275.


\bibitem{Schmidt} Kerz, M; Schmidt, A. Covering Data and Higher Dimensional
Global Class Field Theory. arXive: math/0804.3419.

\bibitem{Lang} Lang, S. Unramified Class Field Theory Over Function Fields
in Several Variables. Annals of Math, 2nd Ser., Vol 64, No. 2 (1956),
285-325.


\bibitem{GIT} Mumford, D; Fogarty, J; Kirwan, F. Geometric Invariant Theory.
Third Enlarged Ed. Springer, Berlin, 1994.

\bibitem{Raskind} Raskind, W. Abelian Calss Field Theory of Arithmetic
Schemes. K-theory and Algebraic Geometry, Proceedings of Symposia in Pure
Mathematics, Vol 58, Part 1 (1995), 85-187.

\bibitem{Saito} Saito, S. Unramified Class Field Theory of Arithmetical
Schemes. Annals of Math, 2nd Ser., Vol 121, No. 2 (1985), 251-281.

\bibitem{VS1} Suslin, A; Voevodsky, V. Singular homology of abstract
algebraic varieties. Invent. Math. 123 (1996), 61-94.

\bibitem{SV2} Suslin, A; Voevodsky, V. Relative Cycles and Chow Sheaves, in
\emph{Cycles, Transfers, and Motivic Homology Theories}, Voevodsky, V;
Suslin, A; Friedlander, E M. Annals of Math Studies, Vol 143. Princeton
University Press, Princeton, NJ, 2000.

\bibitem{w1} Wiesend, G. A Construction of Covers of Arithmetic Schemes. J.
Number Theory, Vol 121 (2006), No. 1, 118-131.

\bibitem{w2} Wiesend, G. Class Field Theory for Arithmetic Schemes. Math
Zeit, Vol 256 (2007), No. 4, 717-729.
\end{thebibliography}
\end{document}